\documentclass[10pt]{article}
\usepackage{amsmath}
\usepackage{amsthm}
\usepackage{amsfonts}
\usepackage{tikz}
\usepackage{mathrsfs}
\usepackage{amsfonts}
\usepackage{amssymb,amsfonts,amsmath, psfrag,eepic,colordvi,epsfig}
 
\topskip 
\numberwithin{equation}{section}
\newtheorem{theorem}{Theorem}[section]

\newcommand{\Rmnum}[1]{\uppercase\expandafter{\romannumeral #1}}

\begin{document}
	\parskip 7pt
	
	\pagenumbering{arabic}
	\def\sof{\hfill\rule{2mm}{2mm}}
	\def\ls{\leq}
	\def\gs{\geq}
	\def\SS{\mathcal S}
	\def\qq{{\bold q}}
	\def\MM{\mathcal M}
	\def\TT{\mathcal T}
	\def\EE{\mathcal E}
	\def\lsp{\mbox{lisp}}
	\def\rsp{\mbox{rasp}}
	\def\pf{\noindent {\it Proof.} }
	\def\mp{\mbox{pyramid}}
	\def\mb{\mbox{block}}
	\def\mc{\mbox{cross}}
	\def\qed{\hfill \rule{4pt}{7pt}}
	\def\pf{\noindent {\it Proof.} }
	\textheight=22cm

	\begin{center}
		{\Large\bf  Combinatorial proof of an inequality on some 
		 partitions separated by parity}
	\end{center}
	
	\begin{center}

	Yan Fan$^{1}$ and Ernest X.W. Xia$^{2}$
	
	$^{1}$School of Mathematical Sciences, \\
	 Jiangsu University, \\
	 Zhenjiang 212013, Jiangsu, People’s Republic
	of China
		
		$^{2}$School of Mathematical Sciences, \\
		Suzhou University of Science and
		Technology, \\
		Suzhou,  215009, Jiangsu Province,
		People’s Republic
		of China

		email:  free1002@126.com,    ernestxwxia@163.com
		
	\end{center}

	\noindent {\bf Abstract.}
	 In 2019, Andrews investigated  integer partitions in which all parts of a given
	parity are smaller than those of the opposite parity and  introduced  
	eight partition functions  based on the parity of the smaller parts and parts of a given
	parity appearing at most once or an unlimited number of times. 
	 Recently, Bringmann, Craig and Nazaroglu
	   studied the asymptotic behavior of the eight partition functions  proved  several inequalities  for sufficiently large $n$. At the end of their paper,
	 they   asked for combinatorial 
	  proofs of those inequalities.  
	   In this paper, we prove that  an inequality   on 
	    partitions separated by parity  holds for $n\geq 373$ by a    combinatorial method. This  answers a question
	     posed  by Bringmann, Craig and Nazaroglu.



	\section{Introduction}
	\allowdisplaybreaks

\raggedbottom

The aim of this paper is to 
 present a combinatorial proof of an inequality  
  involving the number of partitions with parts separated by parity.
   Let us begin with some definitions.   A partition of a positive
 integer
 $n$ is an integer sequence $\lambda
 =(\lambda_1, \lambda_2, \ldots, \lambda_k)$ such that
 $\lambda_1\geq \lambda_2\geq \cdots
 \geq \lambda_k>0$
 and $\lambda_1
 +\lambda_2+\cdots +\lambda_k=n$.
 The numbers $\lambda_i$ are called
 parts of $\lambda$.  The frequency of a part $\lambda_i
 $, denoted
 by $f_{\lambda_i}(\lambda)$,
 is the number of its
 occurrence in the partition $\lambda$ 
  \cite{Andrews-1976}.
 
 In 2019, Andrews \cite{Andrews-2019}  considered partitions in which parts of a given parity are
 all smaller than those of the other parity.
Following the notation in \cite{Andrews-2019}, 
 we denote by $p_{yz}^{wx}(n)$ the number of partitions
    of $n$ with parts of type $w$ satisfying condition $x$ or of type $y$ satisfying
condition $z$ and all parts of type $w$ satisfying condition $x $ are larger than all parts
of type $y$ satisfying condition $z$. Here,  the
  symbols $x$ and $z$ are $u$ or $d$
   (which represent unrestricted parts or distinct parts). 
    The symbols $y$ and $w$ are $e$ or $o$ (which represent even or odd parts)  with $y\neq w$.
  In this paper, we allow partitions counted by $p_{yz}^{wx}(n)
   $ to have no parts of type $y$ or no parts of type $w$.
    For 
example, $p_{od}^{eu}(n)$ denotes 
  the number of partitions of $n$ in which all even parts are greater than all odd parts, even parts are unrestricted and odd parts are distinct. For instance,  $p_{od}^{eu}(5)=3$
   and  the partitions in question  are 
   \[
   (5),\quad (4,1),\quad (2,2,1).
   \]
Andrews \cite{Andrews-2019} considered   all eight families of
partitions, namely $p_{ed}^{od}(n)$, 
$p_{od}^{ed}(n)$,
$p_{od}^{eu}(n)$,
$p_{eu}^{od}(n)$,
$p_{ed}^{ou}(n)$,
$p_{eu}^{ou}(n)$,
$p_{ou}^{ed}(n)$ and $p_{ou}^{eu}(n)$. He also found that 
 some of the eight functions are  closely related to
  mock theta functions. Bringmann and Jennings-Shaffer \cite{Bringmann-2} established 
   new forms of the generating 
    functions for $p_{ed}^{od}(n)$,
     $p_{od}^{ed}(n)$ and  $p_{ou}^{ed}(n)$.
      Recently, Bringmann, Craig and Nazaroglu
   \cite{Bringmann} used Ingham’s Tauberian theorem to compute the asymptotic main
   term for each of the eight functions studied by Andrews
    \cite{Andrews-2019}. Based on the  asymptotic formulas 
     of the eight functions, they  proved that 
    for sufficiently large $n$, 
\[
p_{ed}^{od}(n)<p_{od}^{ed}(n)<
p_{od}^{eu}(n)<
p_{eu}^{od}(n)<
p_{ed}^{ou}(n)<
p_{eu}^{ou}(n)<
p_{ou}^{ed}(n)<p_{ou}^{eu}(n).\]
Note that the second, fifth, and last
 inequalities are 
trivial. At the end their paper, 
Bringmann, Craig and Nazaroglu \cite{Bringmann} asked for combinatorial 
 proofs of the other inequalities.
  Motivated  by Bringmann, Craig and Nazaroglu's open problem,
   Ballantine and Welch \cite{Ballantine} gave combinatorial 
    injections to prove  five inequalities. More precisely, 
     they proved the following theorem.
     
 \begin{theorem} \cite[Theorem 1.1]{Ballantine}
 The following inequalities hold:
 \begin{align*}
 	p_{ed}^{od}(n)&<p_{od}^{ed}(n), \quad {\rm for}\ n\geq 11,\\[6pt]
 	p_{eu}^{ou}(n)&<p_{ou}^{eu}(n), \quad {\rm for}\ n\geq 3,\\[6pt]
 	p_{od}^{eu}(n)&<p_{ed}^{ou}(n), \quad {\rm for}\ n\geq 5,\\[6pt]
 	p_{eu}^{od}(n)&<p_{ou}^{ed}(n), \quad {\rm for}\ n\geq 2,\\[6pt]
 	p_{od}^{ed}(n)&<p_{eu}^{od}(n), \quad {\rm for}\ n\geq 8.
 \end{align*}
 \end{theorem}
 
 Very recently, Fu and Tang \cite{Fu} considered
 partitions with parts separated by parity with
 additional conditions on the multiplicity of parts and posed a conjecture which was confirmed 
  by Ballantine and Welch \cite{Ballantine}.
   At the end of their paper,  Ballantine and Welch \cite{Ballantine} asserted that they were
    unable to find combinatorial arguments for the following inequalities:
    \begin{align}
    p_{od}^{eu}(n)&<p_{eu}^{od}(n), \label{1-1}\\
    p_{eu}^{od}(n)&<p_{ed}^{ou}(n),\label{1-2}\\
      p_{eu}^{ou}(n)&<p_{ou}^{ed}(n).\label{1-3}
    \end{align}
  
  In this paper, we show that  \eqref{1-1} is true 
   for   $n\geq 373$ by a combinatorial method.

 \begin{theorem} \label{Th-1}  
 	For $n\geq 373$, \eqref{1-1} is true.
  \end{theorem}
  
From \cite{Bringmann-2}, we find that   the generating functions of $p_{eu}^{od}(n)$
 and $p_{od}^{eu}(n)$ are  
   \begin{align*}
   \sum_{n=0}^\infty p_{eu}^{od}(n) q^n&=\frac{1}{(q^2;q^2)_\infty}
   \sum_{n=0}^\infty q^{n^2}
   ,\\
  \sum_{n=0}^\infty p_{od}^{eu}(n) (-1)^n 
   q^n&=\frac{1}{(q^2;q^2)_\infty}
  -\frac{1}{(q^2;q^2)_\infty}\sum_{n=1}^\infty\sum_{j=1}^n
   (-1)^{n+j}(1-q^{2n+1})q^{n(3n+1)/2-j^2},
  \end{align*}
  where 
  \[
  (q;q)_\infty:=\prod_{n=1}^\infty (1-q^n).
  \]
With Maple, we find that 
\begin{align*}
	\sum_{n=0}^\infty (p_{eu}^{od}(n)-p_{od}^{eu}(n)) q^n
&=-q^3-q^5-2q^7-q^8-2q^9-4q^{11}-q^{12}-4q^{13}-8q^{15}
-8q^{17}\nonumber\\
&\quad +2q^{18}-\cdots-2q^{49}+816q^{50}+18q^{51}\cdots
\end{align*}  
and for $50\leq n \leq 400$,
  \begin{align}\label{1-4}
  p_{eu}^{od}(n)>p_{od}^{eu}(n).
  \end{align}
 Thanks to Theorem \ref{Th-1} and \eqref{1-4}, we deduce that 
   for $n\geq 50$, 
  \[
  p_{od}^{eu}(n) < p_{eu}^{od}(n).
  \]

  The rest of this paper is organized as
 follows. In Section  2,  we present some definitions on the theory 
  of integer partitions which will be used to prove Theorem \ref{Th-1}.
   Section 3 is devoted to  posing  a combinatorial proof of Theorem \ref{Th-1}.
             We conclude in
      the final section with   a summary of our findings and   a discussion
      of further studies.

 \section{Definitions and notation}
   
     For  a partition $\lambda$, we use the following
      notation:
     
\begin{enumerate}
    
    \item[$\bullet$]
    $\ell(\lambda)$ is  the length of a partition $\lambda$, namely, $\ell(\lambda)$
     is the number of parts in $\lambda$;
     
      \item[$\bullet$]
     $\lambda^{o}$ (resp. $\lambda^{e}$) is the partition whose parts are precisely the odd
      (resp. even) parts of $\lambda$;
      
      \item[$\bullet$]
      $\ell_{o} (\lambda) $ (resp. $\ell_{e} (\lambda) $) is the number of odd (resp. even) parts in $\lambda$;

        \item[$\bullet$] 
        $\lambda_{1}^{o}$ (resp. $\lambda_{1}^{e}$) is  the first part in $\lambda^{o}$
        (resp. $\lambda^{e}$);
        
           \item[$\bullet$] $f_j(\lambda)$ is
             the frequency of the part $j$ in the partition $\lambda$, i.e.,   the number of times $j$  occurs as a
             part in $\lambda$;

                 \end{enumerate}
                 
                   Throughout  this paper, let $A_{od}^{eu}(n)$ (resp. $B_{eu}^{od}(n)$)
                 be the set of partitions of $n$ counted by $p_{od}^{eu}(n)$ (resp. $p_{eu}^{od}(n)$).
                 For a partition
                 $\lambda\in A_{od}^{eu}(n)$, we write it as
                 \[
                 \lambda=(\lambda_1,\lambda_2,\ldots,\lambda_{\ell(\lambda)})=
                 (\lambda_1^{e},\lambda_2^{e},
                 \ldots, \lambda_s^{e},\lambda_1^{o},
                 \ldots, \lambda_t^{o}),
                 \]
                 where $s=\ell_{e} (\lambda)$ and  $t=\ell_{o} (\lambda)$.
                 Note that  $\lambda^{e}=
                 \emptyset$ (resp. $\lambda^{o} =\emptyset$) if and only if $s=0$ (resp. $t=0$). 
                 For a partition
                 $\mu\in B_{eu}^{od}(n)$, we write it as
                 \[
                 \mu=(\mu_1,\mu_2,\ldots,\mu_{\ell(\mu)})=
                 (\mu_1^{o},\mu_2^{o},
                 \ldots, \mu_v^{o},\mu_1^{e},
                 \ldots, \mu_u^{e}),
                 \]
                 where  $u=\ell_{e} (\mu)$ and $v=\ell_{o} (\mu)$.
                 Similarly,  $\mu^{e}=
                 \emptyset$ (resp. $\mu^{o} =\emptyset$) if and only if $u=0$ (resp. $v=0$).

        Occasionally, we use the frequency
        notation for partitions and write $\lambda=(i_1^{f_1},i_2^{f_2},\ldots,i_k^{f_k})$,
         where $f_j:=f_{i_j}(\lambda)$. If $\alpha=(\alpha_1,\alpha_2,\ldots, \alpha_k)$
          and $\beta=(\beta_1,\beta_2,\ldots, \beta_m)$ are two partitions
           with $\alpha_k\geq \beta_1$, then 
           \[
           \alpha\cup \beta=(\alpha_1, \alpha_2, \ldots, \alpha_k,\beta_1,\beta_2,\ldots,\beta_m).
           \]
           
         The Ferrers diagram of a partition $\alpha=(\alpha_1,\alpha_2,\ldots, \alpha_k)$ is
        an array of left justified unit squares
        such that the $i$th row from the top contains $\alpha_i$ unit squares.  
         For example, the Ferrers diagram of $\alpha=(6,4,3,3,1)$ is shown below
      	
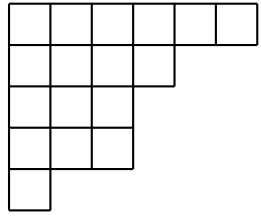
\begin{figure}[htbp] \centering
        	\begin{tikzpicture}[scale=1.1]
        		\draw[thick, black] (8, 0) -- (9.5, 0);
        		\draw[thick, black] (8, 0.5) -- (10, 0.5);
        		\draw[thick, black] (8, 1) -- (11, 1);
        		\draw[thick, black] (8, 1.5) -- (11, 1.5);
        		\draw[thick, black] (8, -0.5) -- (9.5, -0.5);
        		\draw[thick, black] (8, -1) -- (8.5, -1);
        		 
        		\draw[thick, black] (8, -1) -- (8, 1.5);
        		\draw[thick, black] (8.5, -1) -- (8.5, 1.5);
        		\draw[thick, black] (9, -0.5) -- (9, 1.5);
        		\draw[thick, black] (9.5, -0.5) -- (9.5, 1.5);
        		\draw[thick, black] (10, 0.5) -- (10, 1.5);
        		\draw[thick, black] (10.5, 1) -- (10.5, 1.5);
        		\draw[thick, black] (11, 1) -- (11, 1.5);
        		 
        		        	\end{tikzpicture}
        \end{figure}

        \section{Combinatorial proof of Theorem \ref{Th-1}}

  For a partition $\lambda \in A_{od}^{eu}(n)$, we will define an injection
 \[
\psi: A_{od}^{eu}(n)\rightarrow B_{eu}^{od}(n) 
\]
such that $\mu=\psi(\lambda)$. 
  We will break the definition into 17 cases. 
   For $1\leq j \leq 17$,  define
   \begin{align*}
   A_j(n):&=\{\lambda| \lambda  \ {\rm is \ a \ partition
   	 \ of} \ n\ {\rm and \ statifies\  the \ 
   	conditions\  given\  in\  Case } \ j\},\qquad \\
   	B_j(n):&=\{\mu=\psi(\lambda)| \lambda  \in A_j(n)\}.
   \end{align*}
    Due to the discussion of 17 cases  in the proof, we have made  the following  table
      to facilitate readers' understanding of the proof process.
      
 \begin{tabular}{|c|p{4.5cm}|p{8cm}|} 
 	\hline
 	\  &  properties of  $\lambda$ which is  in $A_j(n)$ & properties
 	 of   $\mu=\psi(\lambda) $ which is  in $B_j(n)$\\
 	\hline
 	
 	Case 1  & $s =0 $  or  $t =0$ & $u=0 $  or $  v =0$  \\
 	\hline
 	
 	Case 2  &  $s =t\geq 2 $
 	& 
 	   $u= v\geq 2,  \ \mu_{v}^{o}-\mu_{1}^{e} \geq 2v-3$
  	\\
 	\hline

Case 3 &  $s > t\geq 2$
& 
  $u>v \geq 2,  \mu_1^{o}-\mu_{2}^{o} \geq 2(u-v+1),
  \ \mu_{u-v}^{e}-\mu_{u-v+1}^{e}\geq 
  2v-4 $
\\
\hline

 	 	Case 4  & $	t>s\geq 2 $
 	& 
 	 $v>u \geq 2,
 	 \mu_1^{o}-\mu_{2}^{o} \geq 2(v-u+1),
 	 \ \mu_{v}^{o}-\mu_{1}^{e}\geq 
 	 2u-3 $
 	\\
 	\hline

 	Case 5 & 
 $s=1$, $t\geq 1$
 and  $\lambda_{1}^{e}-\lambda_{1}^{o}\geq 3$ 
 	& $ u= 1,\  v\geq 1 $
 	\\
 	\hline
 	
 	Case 6 & 
 	$s=1$, $ t\geq 5$ and $ \lambda_{1}^{e}-\lambda_{1}^{o}=1$ &  $v\geq 3, \ 2u-3=\mu_1^{o}, \ u
 	-v\geq 3,  \mu_{1}^{e}-\mu_{2}^{e}\geq2,
 	\mu_{3}^{e}=2 $\\
 	\hline
 	
 	Case 7 & 
 		$s =1$, $ t=3 \ or\ 4 $  and $\lambda_{1}^{e}-\lambda_{1}^{o} =1$ 
 	&   $v=3\ {\rm or} \ 4,\  \mu_{v}^{o}=3,\ \mu_{1}^{e}=2,\  
 	u+2v+1 \geq \mu_1^{o} \geq 2v+1,  u\geq 6, \ 2|u$\\
 	\hline
 	
 	Case 8 & 
 	$s=1,  t=2$ and $\lambda_{1}^{e}-\lambda_{1}^{o} =1$ 
 	&  $v=2, u \geq  4, \mu_1^{o}-\mu_2^{o}=2,
 	\mu_1^{e}=2, \mu_2^{o}-	2u+11>0,\ \mu_1^{o} \geq 5$
 	\\
 	\hline
 	
 	Case 9 & 
 	 $s=t=1$, $\lambda_1^{e}-\lambda_1^{o}=1$
 	&  $\mu= (4k-15,5,3,2,2,2)$  \\
 	\hline

 	Case 10 & 
 	$s=2$ and $t=1$
 	&  $v=3, u\geq 5, \mu_{3}^{o}=5,    \mu_1^e=2,\ 2u+25\geq \mu_1^{o}$ \\
 	\hline

 	Case 11 & 
 	$s\geq 3, \ t=1$ and 
 	$   \lambda_1^o=1$ 	
 	&  $u\geq 2,\ v=1$ \\
 	\hline

 	Case 12 & 
 	$s=3,\ t =1$ and $ \lambda_1^{o}\geq 3$
 	&  $ v=3,\  \mu_3^o\geq 7, \ u\geq 4,\  \mu_{1}^e =2$,\ $ 2u+23\geq \mu_1^{o}$ \\
 	\hline

 	Case 13 & 
 	$s=4, \ t =1$ and $\lambda_1^{o}\geq 3$
 	&  $ u\geq 6, \ v=5, \ \mu_5^o=3, \  \mu_{1}^e=2$,
 	 $2u+27\geq \mu_1^{o}$ \\
 	\hline

 	Case 14 & 
 	$s=5,\ t =1$ and $ \lambda_1^{o}\geq 3$
 	&   $u\geq 6,\  
 	v=5,\  \mu_5^o\geq 5, \  \mu_{1}^e =2, 2u+35\geq \mu_{1}^{o}$
 	\\
 	\hline

 	Case 15 & 
 	$6\leq s\leq 10,\ t =1$ and $ \lambda_1^{o}\geq 3$ & $f_2(\mu)>12,\ v=3,   \  	
 	\mu_{u-f_2(\mu)}^{e}\geq 4, \  3 \leq 	u-f_2(\mu)\leq  7, 
 	2f_2(\mu)+15\geq \mu_1^{o}$ \\
 	\hline

 	Case 16 & 
 	$s\geq 11,\ t=1,\  \lambda_1^{o}\geq 3$ and $ \lambda_1^{e}-\lambda_2^{e}\leq 10$
 	& $u\geq 9,\ 
 	v=3, \  f_2(\mu)\leq 5, \ 	   	\mu_1^{o}-\mu_2^{o}\leq 12
 	,\ \mu_{u-5}^{e}-\mu_{u-4}^{e}\geq 2$  \\
 	\hline

 	Case 17 & 
 	$ s\geq11, \ t=1,\ \lambda_1^{o}\geq 3$ and $\lambda_1^{e}-\lambda_2^{e}\geq 12$
 	& $u\geq 15,\ v=3, \  6 \leq  f_2(\mu) \leq 11, \mu_{u-11}^{e}-\mu_{u-10}^{e}\geq 2$  \\
 	\hline
 	
 \end{tabular}

 From the above table, we know that 
 \begin{align}
 A_{od}^{eu}(n)=\mathop{\bigcup}
 \limits_{j=1}^{17}
 A_j(n).\label{3-1}
 \end{align}
 and  if $i\neq j$, then
 \begin{align}\label{3-2}
  A_i(n)\bigcap A_j(n)=B_i(n)\bigcap B_j(n)=\emptyset.
 \end{align}
 
Note that if $n=6k+j\ (j\in\{1,3,5\})$ with $n\geq 373$, then
 $(2k+3,2k+1,2k+j-8,2,2)\in B_{eu}^{od}(n)$  and $(2k+3,2k+1,2k+j-8,2,2)\not\in \mathop{\bigcup}
 \limits_{j=1}^{17}
 B_j (n)$.  In addition,  if $n=6k+j,\ (j\in\{0,2,4\})$ with $n\geq 373$, then
 $(2k+1,2k-1,2k+j-7,3,2,2)\in B_{eu}^{od}(n)$  and $(2k+1,2k-1,2k+j-7,3,2,2)\not\in \mathop{\bigcup}
 \limits_{j=1}^{17}
 B_j (n)$.
  Therefore, for $n\geq 373$, 
  \begin{align}\label{3-3}
  \mathop{\bigcup}
 \limits_{j=1}^{17}
 B_j (n)\subset B_{eu}^{od}(n) \qquad  {\rm and } 
 \quad \mathop{\bigcup}
 \limits_{j=1}^{17}B_j(n) \neq  B_{eu}^{od}(n).
 \end{align}
In view of \eqref{3-1}--\eqref{3-3},
\begin{align}
 p_{od}^{eu}(n)=|A_{od}^{eu}(n)|=\sum_{j=1}^{17} |A_j(n)|.\label{3-4}
\end{align}
and 
 \begin{align}
 	p_{eu}^{od}(n)=|B_{eu}^{od}(n)|>\sum_{j=1}^{17} |B_j(n)|.\label{3-5}
 \end{align}
  
 Now, we describe   the details of  $\psi$.
 
Case 1:
	 $s =0 $  or  $t=0$. For any $\lambda\in A_1(n)$
	  with $n\geq 1$, define 
	\begin{align*}
		\psi(\lambda):=\lambda
	\end{align*}
	then
	\begin{align*}
		B_{1}(n)=\{\mu \in B_{eu}^{od}(n) \ | u =0  \ {\rm or}\ 
		   v =0 \}
	\end{align*}
	Obviously, for any $\mu\in B_1(n)$ with $n\geq 1$, 
	\begin{align*}
		\varphi_{1}^{-1}(\mu)=\mu.
	\end{align*}

Case 2:
	  $s=t\geq 2$.  For  any $\lambda=(\lambda_1^{e},\ldots,
	  \lambda_s^{e},\lambda_1^{o},\ldots,
	  \lambda_s^{o})\in A_2(n)$
	  with $n\geq 12$,  define 
	\begin{align*}
		\psi(\lambda):=
		 (&\lambda_{1}^{e} +(2s-3), \lambda_{2}^{e}
		 +(2s-5),\ldots, \lambda_{s-1}^{e}+1 ,\lambda_{s}^{e} -1,\nonumber\\
		 &\lambda_{1}^{o}-(2s-3) ,\lambda_{2}^{o}-(2s-5),\ldots,\lambda_{s-1}^{o} -1, \lambda_{s}^{o}+1 ).
		 \end{align*}
Then 
	\begin{align*}
		B_{2}(n)=& \{\mu \in B_{eu}^{od}(n) \ | u= v\geq 2,  \mu_{v}^{o}-\mu_{1}^{e} \geq 2v-3 \}	.
	\end{align*}
	Obviously, for  any $\mu=(\mu_1^{o},\ldots,
	\lambda_u^{o},\mu_1^{e},\ldots,
	\lambda_u^{e})\in B_2(n)$
	with $n\geq 12$,  
	\begin{align*}
		\psi^{-1}(\mu):=
		(&\mu_{1}^{o} -(2u-3), \mu_{2}^{o}
		-(2u-5),\ldots, \mu_{u-1}^{o}-1 ,\mu_{u}^{o} +1,\nonumber\\
		&\mu_{1}^{e}+(2u-3) ,\mu_{2}^{e}+(2u-5),\ldots,\mu_{u-1}^{e} +1, \mu_{u}^{e}-1 ).
	\end{align*}
	Example: \ Let $n=39$ and $\lambda=(8,8,8,7,5,3)$. Then $s=t=3$ and  $\psi(\lambda)=(8+3,8+1,8-1,7-3,5-1,3+1)=(11,9,7,4,4,4).$
	
	
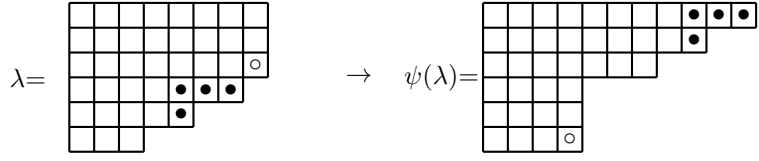
\begin{figure}[htbp] \centering
	 \begin{tikzpicture}[scale=1.1]
			 \draw[thick, black] (4, 0) -- (6.4, 0);
			 \draw[thick, black] (4, 0.3) -- (6.4, 0.3);
			 \draw[thick, black] (4, 0.6) -- (6.4, 0.6);
			 \draw[thick, black] (4, 0.9) -- (6.4, 0.9);
			 \draw[thick, black] (4, -0.3) -- (6.1, -0.3);
			    \draw[thick, black] (4, -0.6) -- (5.5, -0.6);
			   \draw[thick, black] (4, -0.9) -- (4.9, -0.9);
			   
		    \draw[thick, black] (4, -0.9) -- (4, 0.9);
		    \draw[thick, black] (4.3, -0.9) -- (4.3, 0.9);
		     \draw[thick, black] (4.6, -0.9) -- (4.6, 0.9);
		      \draw[thick, black] (4.9, -0.9) -- (4.9, 0.9);
		       \draw[thick, black] (5.2, -0.6) -- (5.2, 0.9);
		        \draw[thick, black] (5.5, -0.6) -- (5.5, 0.9);
		        \draw[thick, black] (5.8, -0.3) -- (5.8, 0.9);
		        \draw[thick, black] (6.1, -0.3) -- (6.1, 0.9);
		        \draw[thick, black] (6.4, 0) -- (6.4, 0.9);
			
			\draw (3.5, 0) node{$\lambda$=};
			\draw (5.35, -0.15) node{$\bullet$};
			\draw (5.65, -0.15) node{$\bullet$};
			\draw (5.95, -0.15) node{$\bullet$};
			\draw (5.35, -0.45) node{$\bullet$};
			\draw (6.25, 0.15) node{$\circ$};

			\draw[thick, black] (9, 0) -- (11.1, 0);
			\draw[thick, black] (9, 0.3) -- (11.7, 0.3);
			\draw[thick, black] (9, 0.6) -- (12.3, 0.6);
			\draw[thick, black] (9, 0.9) -- (12.3, 0.9);
			\draw[thick, black] (9, -0.3) -- (10.2, -0.3);
			\draw[thick, black] (9, -0.6) -- (10.2, -0.6);
			\draw[thick, black] (9, -0.9) -- (10.2, -0.9);
			
			\draw[thick, black] (9, -0.9) -- (9, 0.9);
			\draw[thick, black] (9.3, -0.9) -- (9.3, 0.9);
			
			\draw[thick, black] (9.6, -0.9) -- (9.6, 0.9);
			\draw[thick, black] (9.9, -0.9) -- (9.9, 0.9);
			\draw[thick, black] (10.2, -0.9) -- (10.2, 0.9);
			\draw[thick, black] (10.5, 0) -- (10.5, 0.9);
			\draw[thick, black] (10.8, 0) -- (10.8, 0.9);
			\draw[thick, black] (11.1, 0) -- (11.1, 0.9);
			\draw[thick, black] (11.4, 0.3) -- (11.4, 0.9);
			\draw[thick, black] (11.7, 0.3) -- (11.7, 0.9);
			\draw[thick, black] (12, 0.6) -- (12, 0.9);
			\draw[thick, black] (12.3, 0.6) -- (12.3, 0.9);
			
			\draw (7.5, 0) node{$\rightarrow$};
			\draw (8.5, 0) node{$\psi(\lambda)$=};
			\draw (11.55, 0.75) node{$\bullet$};
			\draw (11.85, 0.75) node{$\bullet$};
			\draw (12.15, 0.75) node{$\bullet$};
			\draw (11.55, 0.45) node{$\bullet$};
			\draw (10.05, -0.75) node{$\circ$};
			 \end{tikzpicture}
			\end{figure}

	\newpage
	Case 3:
	$s>t\geq 2$.  For any $\lambda\in A_3(n)$
	with $n\geq 16$,  define 
	\begin{align*}
		\psi(\lambda):=
		(&\lambda_{1}^{e}  +2(s-t)+(2t-3), \lambda_{2}^{e}
		+(2t-5),\ldots, \lambda_{t-1}^{e}+1 ,\lambda_{t}^{e} -1,
		 \lambda_{t+1}^{e}-2,\lambda_{t+2}^{e}-2,\ldots, \lambda_{s}^{e}-2,\nonumber\\
		&\qquad \qquad \ \lambda_{1}^{o}-(2t-3) ,\lambda_{2}^{o}-(2t-5),\ldots,\lambda_{t-1}^{o} -1, \lambda_{t}^{o}+1 ).
	\end{align*}
	Then 
	\begin{align*}
		B_{3}(n)=& \{\mu \in B_{eu}^{od} \ | u>v \geq 2,  \mu_1^{o}-\mu_{2}^{o} \geq 2(u-v+1),
		\ \mu_{u-v}^{e}-\mu_{u-v+1}^{e}\geq 
		 2v-4 \}	.
	\end{align*}
	Obviously, for any $\mu \in B_3(n)$ with $n\geq 16$, 
	\begin{align*}
		\psi^{-1}(\mu):=
		(&\mu_{1}^{o} -(2(u-v)+2v-3), \mu_{2}^{o}
		-(2v-5),\ldots, \mu_{v-1}^{o}-1 ,\mu_{v}^{o} +1,\mu_{1}^{e}+2,\mu_{2}^{e}+2,
		\ldots, \mu_{u-v}^{e}+2\nonumber\\
		&\mu_{u-v+1}^{e}+(2v-3) ,\mu_{u-v+2}^{e}+(2v-5),\ldots,\mu_{u-1}^{e} +1, \mu_{u}^{e}-1 ).
	\end{align*}

Example. Let $n=26$ and $\lambda=(6,6,6,4,3,1)$. Then $\varphi_{2}(\lambda)=(6+4+1,6-1,6-2,4-2,3-1,1+1)=(11,5,4,2,2,2).$

	
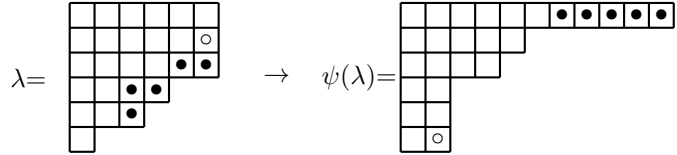
\begin{figure}[htbp] \centering
		\begin{tikzpicture}[scale=1.1]
			\draw[thick, black] (4, 0) -- (5.8, 0);
			\draw[thick, black] (4, 0.3) -- (5.8, 0.3);
			\draw[thick, black] (4, 0.6) -- (5.8, 0.6);
			\draw[thick, black] (4, 0.9) -- (5.8, 0.9);
			\draw[thick, black] (4, -0.3) -- (5.2, -0.3);
			\draw[thick, black] (4, -0.6) -- (4.9, -0.6);
			\draw[thick, black] (4, -0.9) -- (4.3, -0.9);
			
			\draw[thick, black] (4, -0.9) -- (4, 0.9);
			\draw[thick, black] (4.3, -0.9) -- (4.3, 0.9);
			\draw[thick, black] (4.6, -0.6) -- (4.6, 0.9);
			\draw[thick, black] (4.9, -0.6) -- (4.9, 0.9);
			\draw[thick, black] (5.2, -0.3) -- (5.2, 0.9);
			\draw[thick, black] (5.5, 0) -- (5.5, 0.9);
			\draw[thick, black] (5.8, 0) -- (5.8, 0.9);

			\draw (3.5, 0) node{$\lambda$=};
			\draw (4.75, -0.45) node{$\bullet$};
			\draw (4.75, -0.15) node{$\bullet$};
			\draw (5.05, -0.15) node{$\bullet$};
			\draw (5.35, 0.15) node{$\bullet$};
			\draw (5.65, 0.15) node{$\bullet$};
			\draw (5.65, 0.45) node{$\circ$};

			\draw[thick, black] (8, 0) -- (9.2, 0);
			\draw[thick, black] (8, 0.3) -- (9.5, 0.3);
			\draw[thick, black] (8, 0.6) -- (11.3, 0.6);
			\draw[thick, black] (8, 0.9) -- (11.3, 0.9);
			\draw[thick, black] (8, -0.3) -- (8.6, -0.3);
			\draw[thick, black] (8, -0.6) -- (8.6, -0.6);
			\draw[thick, black] (8, -0.9) -- (8.6, -0.9);
			
			\draw[thick, black] (8, -0.9) -- (8, 0.9);
			\draw[thick, black] (8.3, -0.9) -- (8.3, 0.9);

			\draw[thick, black] (8.6, -0.9) -- (8.6, 0.9);
			\draw[thick, black] (8.9, 0) -- (8.9, 0.9);
			\draw[thick, black] (9.2, 0) -- (9.2, 0.9);
			\draw[thick, black] (9.5, 0.3) -- (9.5, 0.9);
			\draw[thick, black] (9.8, 0.6) -- (9.8, 0.9);
			\draw[thick, black] (10.1, 0.6) -- (10.1, 0.9);
			\draw[thick, black] (10.4, 0.6) -- (10.4, 0.9);
			\draw[thick, black] (10.7, 0.6) -- (10.7, 0.9);
			\draw[thick, black] (11, 0.6) -- (11, 0.9);
			\draw[thick, black] (11.3, 0.6) -- (11.3, 0.9);
			
			\draw (6.5, 0) node{$\rightarrow$};
			\draw (7.5, 0) node{$\psi(\lambda)$=};
			\draw (9.95, 0.75) node{$\bullet$};
			\draw (10.25, 0.75) node{$\bullet$};
			\draw (10.55, 0.75) node{$\bullet$};
			\draw (10.85, 0.75) node{$\bullet$};
			\draw (11.15, 0.75) node{$\bullet$};
			\draw (8.45, -0.75) node{$\circ$};
		\end{tikzpicture}
	\end{figure}

		Case 4:
	$t>s\geq 2$.  For any $\lambda\in A_4(n)$ with $n\geq 21$,  define 
	\begin{align*}
		\psi(\lambda):=
		(&\lambda_{1}^{e}  +2(t-s)+(2s-3), \lambda_{2}^{e}
		+(2s-5),\ldots, \lambda_{s-1}^{e}+1 ,\lambda_{s}^{e} -1,
		\lambda_{1}^{o}-2,\lambda_{2}^{o}-2,\ldots, \lambda_{t-s}^{o}-2,\nonumber\\
		&\qquad \ \lambda_{t-s+1}^{o}-(2s-3) ,\lambda_{t-s+2}^{o}-(2s-5),\ldots,\lambda_{t-1}^{o} -1, \lambda_{t}^{o}+1 ).
	\end{align*}
	Therefore,
	\begin{align*}
		B_{4}(n)=& \{\mu \in B_{eu}^{od}(n) \ | v>u \geq 2,
		  \mu_1^{o}-\mu_{2}^{o} \geq 2(v-u+1),
		\ \mu_{v}^{o}-\mu_{1}^{e}\geq 
		2u-3 \}	.
	\end{align*}
	In addition,
	 for any $\mu 
	 \in B_4(n)$ with $n\geq 21$, 
	\begin{align*}
		\psi^{-1}(\mu):=
		(&\mu_{1}^{o} -(2(v-u)+2u-3), \mu_{2}^{o}
		-(2u-5),\ldots, \mu_{u-1}^{o}-1 ,\mu_{u}^{o} +1,\mu_{u+1}^{o}+2,
		\ldots, \mu_{v}^{o}+2\nonumber\\
		&\qquad \qquad \quad \mu_{1}^{e}+(2u-3) ,\mu_{2}^{e}+(2u-5),\ldots,\mu_{u-1}^{e} +1, \mu_{u}^{e}-1 ).
	\end{align*}

Example. Let $n=32$ and 
	 $\lambda=(8,8,7,5,3,1)$. Then$\psi(\lambda)=(8+4+1,8-1,7-2,5-2,3-1,1+1)=(13,7,5,3,2,2)$.

	 
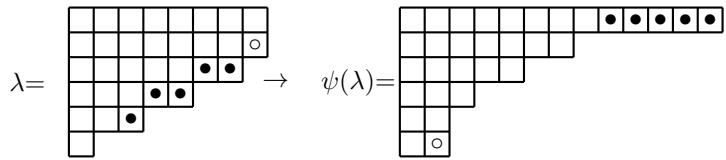
\begin{figure}[htbp]  \centering
	 	\begin{tikzpicture}[scale=1.1]
	 		\draw[thick, black] (4, 0) -- (6.1, 0);
	 		\draw[thick, black] (4, 0.3) -- (6.4, 0.3);
	 		\draw[thick, black] (4, 0.6) -- (6.4, 0.6);
	 		\draw[thick, black] (4, 0.9) -- (6.4, 0.9);
	 		\draw[thick, black] (4, -0.3) -- (5.5, -0.3);
	 		\draw[thick, black] (4, -0.6) -- (4.9, -0.6);
	 		\draw[thick, black] (4, -0.9) -- (4.3, -0.9);
	 		
	 		\draw[thick, black] (4, -0.9) -- (4, 0.9);
	 		\draw[thick, black] (4.3, -0.9) -- (4.3, 0.9);
	 		\draw[thick, black] (4.6, -0.6) -- (4.6, 0.9);
	 		\draw[thick, black] (4.9, -0.6) -- (4.9, 0.9);
	 		\draw[thick, black] (5.2, -0.3) -- (5.2, 0.9);
	 		\draw[thick, black] (5.5, -0.3) -- (5.5, 0.9);
	 		\draw[thick, black] (5.8, 0) -- (5.8, 0.9);
	 		\draw[thick, black] (6.1, 0) -- (6.1, 0.9);
	 		\draw[thick, black] (6.4, 0.3) -- (6.4, 0.9);

	 		\draw (3.5, 0) node{$\lambda$=};
	 		\draw (4.75, -0.45) node{$\bullet$};
	 		\draw (5.05, -0.15) node{$\bullet$};
	 		\draw (5.35, -0.15) node{$\bullet$};
	 		\draw (5.65, 0.15) node{$\bullet$};
	 		\draw (5.95, 0.15) node{$\bullet$};
	 		\draw (6.25, 0.45) node{$\circ$};

	 		\draw[thick, black] (8, 0) -- (9.5, 0);
	 		\draw[thick, black] (8, 0.3) -- (10.1, 0.3);
	 		\draw[thick, black] (8, 0.6) -- (11.9, 0.6);
	 		\draw[thick, black] (8, 0.9) -- (11.9, 0.9);
	 		\draw[thick, black] (8, -0.3) -- (8.9, -0.3);
	 		\draw[thick, black] (8, -0.6) -- (8.6, -0.6);
	 		\draw[thick, black] (8, -0.9) -- (8.6, -0.9);
	 		
	 		\draw[thick, black] (8, -0.9) -- (8, 0.9);
	 		\draw[thick, black] (8.3, -0.9) -- (8.3, 0.9);

	 		\draw[thick, black] (8.6, -0.9) -- (8.6, 0.9);
	 		\draw[thick, black] (8.9, -0.3) -- (8.9, 0.9);
	 		\draw[thick, black] (9.2, 0) -- (9.2, 0.9);
	 		\draw[thick, black] (9.5, 0) -- (9.5, 0.9);
	 		
	 		\draw[thick, black] (9.8, 0.3) -- (9.8, 0.9);
	 		\draw[thick, black] (10.1, 0.3) -- (10.1, 0.9);
	 		\draw[thick, black] (10.4, 0.6) -- (10.4, 0.9);
	 		\draw[thick, black] (10.7, 0.6) -- (10.7, 0.9);
	 		\draw[thick, black] (11, 0.6) -- (11, 0.9);
	 		\draw[thick, black] (11.3, 0.6) -- (11.3, 0.9);
	 		\draw[thick, black] (11.6, 0.6) -- (11.6, 0.9);
	 		\draw[thick, black] (11.9, 0.6) -- (11.9, 0.9);

	 		\draw (6.5, 0) node{$\rightarrow$};
	 		\draw (7.5, 0) node{$\psi(\lambda)$=};
	 		\draw (10.55, 0.75) node{$\bullet$};
	 		\draw (10.85, 0.75) node{$\bullet$};
	 		\draw (11.15, 0.75) node{$\bullet$};
	 		\draw (11.45, 0.75) node{$\bullet$};
	 		\draw (11.75, 0.75) node{$\bullet$};
	 		\draw (8.45, -0.75) node{$\circ$};
	 	\end{tikzpicture}
	 \end{figure}
	 
	 \newpage
	 
Case 5: $s=1$, $t\geq 1$
 and  $\lambda_{1}^{e}-\lambda_{1}^{o}\geq 3.$ 
	 For any $\lambda\in A_5(n)$ with $n\geq 5$, define 
	 	\[
	 		\psi(\lambda):=(\lambda_1^{e}-1,\lambda_1^{o},\ldots,\lambda_{t-1}^{o},
	 		\lambda_t^{o}+1).
	 	\]
Therefore, 
	 	\begin{align*}
	 		B_{5}(n)=\{ \mu \in B_{eu}^{od}(n) \ | u= 1, v\geq 1
	 		\}.
	 	\end{align*}
Furthermore, for any $\mu\in B_5(n)$ with $n\geq 5$,  
	 		\begin{align*}
	 		\psi^{-1}(\mu):=(\mu_1^{o}+1,\mu_2^{o},\ldots,\mu_{v}^{o},
	 		\mu_1^{e}-1).
	 	\end{align*}
	 	
	 	Example.  Let $n=10$ and $\lambda=(6,3,1)$. Then $\psi(\lambda)=(6-1,3,1+1)=(5,3,2).$

Case 6:	$s=1$, $ t\geq 5$ and $ \lambda_{1}^{e}-\lambda_{1}^{o}=1$.  
Note that $\lambda_{1}^{e}\geq 2t$. 
  For any $\lambda\in A_6(n)$ with $n\geq 35$, define  
	 	\begin{align*}
	 		\psi(\lambda):=(\lambda_1^{o},\lambda_2^{o},\ldots,\lambda_{t-2}^{o},
	 		\lambda_{t-1}^{o}+1,\lambda_{t}^{o}+1
	 		)\cup (2^{(\lambda_1^{e}-2)/2}).
	 	\end{align*}
Thus, 
	 	\begin{align*}
	 		B_{6}(n)=\{ & \mu \in B_{eu}^{od}(n) \ | v\geq 3, \ 2u-3=\mu_1^{o}, \ u
	 		-v\geq 3,  \mu_{1}^{e}-\mu_{2}^{e}\geq2,
	 		\mu_{3}^{e}=2
	 		\}.
	 	\end{align*}
  	Furthermore, 
 for any $\mu\in B_6(n)$ with $n\geq 35$, 
	 	\begin{align*}
	\varphi^{-1}(\mu)=\big(2u-2,\mu_1^{o},\mu_2^{o},\ldots,\mu_{v}^{o},
	\mu_{1}^{e}-1,\mu_{2}^{e}-1
	\big).
\end{align*}

	 	Example.  Let $n=35$ and $\lambda=(10,9,7,5,3,1)$. Then $\psi(\lambda)=(9,7,5,3+1,1+1)\cup(2,2,2,2)=(9,7,5,4,2,2,2,2,2).$
	 	
	 	 
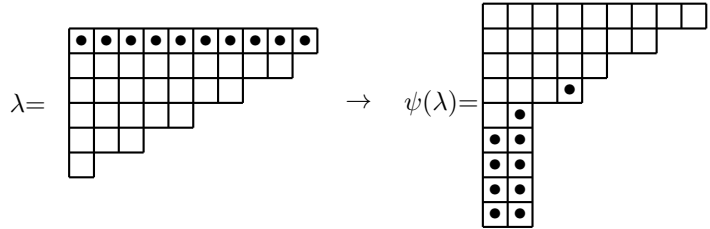
\begin{figure}[htbp] \centering
	 		\begin{tikzpicture}[scale=1.1]
	 			\draw[thick, black] (4, 0) -- (6.1, 0);
	 			\draw[thick, black] (4, 0.3) -- (6.7, 0.3);
	 			\draw[thick, black] (4, 0.6) -- (7, 0.6);
	 			\draw[thick, black] (4, 0.9) -- (7, 0.9);
	 			\draw[thick, black] (4, -0.3) -- (5.5, -0.3);
	 			\draw[thick, black] (4, -0.6) -- (4.9, -0.6);
	 			\draw[thick, black] (4, -0.9) -- (4.3, -0.9);
	 			
	 			\draw[thick, black] (4, -0.9) -- (4, 0.9);
	 			\draw[thick, black] (4.3, -0.9) -- (4.3, 0.9);
	 			\draw[thick, black] (4.6, -0.6) -- (4.6, 0.9);
	 			\draw[thick, black] (4.9, -0.6) -- (4.9, 0.9);
	 			\draw[thick, black] (5.2, -0.3) -- (5.2, 0.9);
	 			\draw[thick, black] (5.5, -0.3) -- (5.5, 0.9);
	 			\draw[thick, black] (5.8, 0) -- (5.8, 0.9);
	 			\draw[thick, black] (6.1, 0) -- (6.1, 0.9);
	 			\draw[thick, black] (6.4, 0.3) -- (6.4, 0.9);
	 			\draw[thick, black] (6.7, 0.3) -- (6.7, 0.9);
	 			\draw[thick, black] (7, 0.6) -- (7, 0.9);
	 			
	 			\draw (3.5, 0) node{$\lambda$=};
	 			\draw (4.15, 0.75) node{$\bullet$};
	 			\draw (4.45, 0.75) node{$\bullet$};
	 			\draw (4.75, 0.75) node{$\bullet$};
	 			\draw (5.05, 0.75) node{$\bullet$};
	 			\draw (5.35, 0.75) node{$\bullet$};
	 			\draw (5.65, 0.75) node{$\bullet$};
	 			\draw (5.95, 0.75) node{$\bullet$};
	 			\draw (6.25, 0.75) node{$\bullet$};
	 			\draw (6.55, 0.75) node{$\bullet$};
	 			\draw (6.85, 0.75) node{$\bullet$};

	 			\draw[thick, black] (9, 0) -- (10.2, 0);
	 			\draw[thick, black] (9, 0.3) -- (10.5, 0.3);
	 			\draw[thick, black] (9, 0.6) -- (11.1, 0.6);
	 			\draw[thick, black] (9, 0.9) -- (11.7, 0.9);
	 			\draw[thick, black] (9, 1.2) -- (11.7, 1.2);
	 			\draw[thick, black] (9, -0.3) -- (9.6, -0.3);
	 			\draw[thick, black] (9, -0.6) -- (9.6, -0.6);
	 			\draw[thick, black] (9, -0.9) -- (9.6, -0.9);
	 			\draw[thick, black] (9, -1.2) -- (9.6, -1.2);
	 			\draw[thick, black] (9, -1.5) -- (9.6, -1.5);
	 			
	 			\draw[thick, black] (9, -1.5) -- (9, 1.2);
	 			\draw[thick, black] (9.3, -1.5) -- (9.3, 1.2);

	 			\draw[thick, black] (9.6, -1.5) -- (9.6, 1.2);
	 			\draw[thick, black] (9.9, 0) -- (9.9, 1.2);
	 			\draw[thick, black] (10.2, 0) -- (10.2, 1.2);
	 			\draw[thick, black] (10.5, 0.3) -- (10.5, 1.2);
	 			\draw[thick, black] (10.8, 0.6) -- (10.8, 1.2);
	 			\draw[thick, black] (11.1, 0.6) -- (11.1, 1.2);
	 			\draw[thick, black] (11.4, 0.9) -- (11.4, 1.2);
	 			\draw[thick, black] (11.7, 0.9) -- (11.7, 1.2); 
	 			
	 			\draw (7.5, 0) node{$\rightarrow$};
	 			\draw (8.5, 0) node{$\psi(\lambda)$=};
	 			\draw (10.05, 0.15) node{$\bullet$};
	 			\draw (9.45, -0.15) node{$\bullet$};
	 			\draw (9.15, -0.45) node{$\bullet$};
	 			\draw (9.45, -0.45) node{$\bullet$};
	 			
	 				\draw (9.15, -0.75) node{$\bullet$};
	 			\draw (9.45, -0.75) node{$\bullet$};
	 			
	 			\draw (9.15, -1.05) node{$\bullet$};
	 		\draw (9.45, -1.05) node{$\bullet$};

	 				\draw (9.15, -1.35) node{$\bullet$};
	 			\draw (9.45, -1.35) node{$\bullet$};
	 			
	 			\end{tikzpicture}
	 	\end{figure}

	Case 7: 		$s =1$, $ t=3 \ or\ 4 $  and $\lambda_{1}^{e}-\lambda_{1}^{o} =1.$ 
	 For $\lambda\in A_7(n) $ with $n \geq 54 $, 
	 			define 
	 		\begin{align*}
	 			\psi(\lambda):=(\lambda_2^{o}+4,\ldots,\lambda_t^{o}+4,3)
	 			\cup (2^{\lambda_1^{e}-2t}).
	 		\end{align*}
	 	 Note that if $n\geq 54$, then $\lambda_1^{e}\geq 14$ and 
	 		\begin{align*}
	 			B_{7}(n)=\{ &\mu \in B_{eu}^{od} (n) \ | v=3\ {\rm or} \ 4,\  \mu_{v}^{o}=3,\ \mu_{1}^{e}=2,\  
	 			u+2v+1 \geq \mu_1^{o} \geq 2v+1,  u\geq 6, \ 2|u 
	 			\}.
	 		\end{align*}
	 	Moreover,
	 	 for any $\mu\in B_7(n)$ with $n\geq 54$,
	 	 \begin{align*}
	 	 	\psi^{-1}(\mu):=(u+2v,u+2v-1, \mu_1^{o}-4,\ldots,\mu_{v-1}^{o}-4).
	 	 \end{align*}	
	 		
	 		Example.  Let $n=54$ and $\lambda=(14,13,11,9,7)$. Then $\psi(\lambda)=(11+4,9+4,7+4,3)\cup (2^{6})=(15,13,11,3,2,2,2,2,2,2).$
	 		 
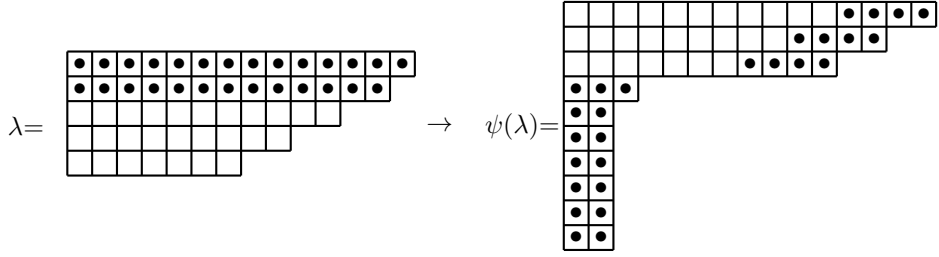
	 \begin{figure}[htbp] \centering
	 	\begin{tikzpicture}[scale=1.1]
	 		\draw[thick, black] (4, 0) -- (7.3, 0);
	 		\draw[thick, black] (4, 0.3) -- (7.9, 0.3);
	 		\draw[thick, black] (4, 0.6) -- (8.2, 0.6);
	 		\draw[thick, black] (4, 0.9) -- (8.2, 0.9);
	 		\draw[thick, black] (4, -0.3) -- (6.7, -0.3);
	 		\draw[thick, black] (4, -0.6) -- (6.1, -0.6);
	 		
	 		\draw[thick, black] (4, -0.6) -- (4, 0.9);
	 		\draw[thick, black] (4.3, -0.6) -- (4.3, 0.9);
	 		\draw[thick, black] (4.6, -0.6) -- (4.6, 0.9);
	 		\draw[thick, black] (4.9, -0.6) -- (4.9, 0.9);
	 		\draw[thick, black] (5.2, -0.6) -- (5.2, 0.9);
	 		\draw[thick, black] (5.5, -0.6) -- (5.5, 0.9);
	 		\draw[thick, black] (5.8, -0.6) -- (5.8, 0.9);
	 		\draw[thick, black] (6.1, -0.6) -- (6.1, 0.9);
	 		\draw[thick, black] (6.4,-0.3) -- (6.4, 0.9);
	 		\draw[thick, black] (6.7, -0.3) -- (6.7, 0.9);
	 		\draw[thick, black] (7, 0) -- (7, 0.9);
	 		
	 		\draw[thick, black] (7.3, 0) -- (7.3, 0.9);
	 		\draw[thick, black] (7.6, 0.3) -- (7.6, 0.9);
	 		\draw[thick, black] (7.9, 0.3) -- (7.9, 0.9);
	 		\draw[thick, black] (8.2, 0.6) -- (8.2, 0.9);
	 		
	 		\draw (3.5, 0) node{$\lambda$=};
	 		\draw (4.15, 0.75) node{$\bullet$};
	 		\draw (4.45, 0.75) node{$\bullet$};
	 		\draw (4.75, 0.75) node{$\bullet$};
	 		\draw (5.05, 0.75) node{$\bullet$};
	 		\draw (5.35, 0.75) node{$\bullet$};
	 		\draw (5.65, 0.75) node{$\bullet$};
	 		\draw (5.95, 0.75) node{$\bullet$};
	 		\draw (6.25, 0.75) node{$\bullet$};
	 		\draw (6.55, 0.75) node{$\bullet$};
	 		\draw (6.85, 0.75) node{$\bullet$};
	 		
	 		\draw (7.15, 0.75) node{$\bullet$};
	 		\draw (7.45, 0.75) node{$\bullet$};
	 		\draw (7.75, 0.75) node{$\bullet$};
	 		\draw (8.05, 0.75) node{$\bullet$};
	 		
	 	\draw (4.15, 0.45) node{$\bullet$};
	 	\draw (4.45, 0.45) node{$\bullet$};
	 	\draw (4.75, 0.45) node{$\bullet$};
	 	\draw (5.05, 0.45) node{$\bullet$};
	 	\draw (5.35, 0.45) node{$\bullet$};
	 	\draw (5.65, 0.45) node{$\bullet$};
	 	\draw (5.95, 0.45) node{$\bullet$};
	 	\draw (6.25, 0.45) node{$\bullet$};
	 	\draw (6.55, 0.45) node{$\bullet$};
	 	\draw (6.85, 0.45) node{$\bullet$};
	 	
	 	\draw (7.15, 0.45) node{$\bullet$};
	 	\draw (7.45, 0.45) node{$\bullet$};
	 	\draw (7.75, 0.45) node{$\bullet$};

	 		\draw[thick, black] (10, 0) -- (10.6, 0);
	 		\draw[thick, black] (10, 0.3) -- (10.9, 0.3);
	 		\draw[thick, black] (10, 0.6) -- (13.3, 0.6);
	 		\draw[thick, black] (10, 0.9) -- (13.9, 0.9);
	 		\draw[thick, black] (10, 1.2) -- (14.5, 1.2);
	 		\draw[thick, black] (10, 1.5) -- (14.5, 1.5);
	 		\draw[thick, black] (10, -0.3) -- (10.6, -0.3);
	 		\draw[thick, black] (10, -0.6) -- (10.6, -0.6);
	 		\draw[thick, black] (10, -0.9) -- (10.6, -0.9);
	 		\draw[thick, black] (10, -1.2) -- (10.6, -1.2);
	 		\draw[thick, black] (10, -1.5) -- (10.6, -1.5);

	 		\draw[thick, black] (10, -1.5) -- (10, 1.5);
	 		\draw[thick, black] (10.3, -1.5) -- (10.3, 1.5);
	 		\draw[thick, black] (10.6, -1.5) -- (10.6, 1.5);
	 		\draw[thick, black] (10.9, 0.3) -- (10.9, 1.5);
	 		\draw[thick, black] (11.2, 0.6) -- (11.2, 1.5);
	 		\draw[thick, black] (11.5, 0.6) -- (11.5, 1.5);
	 		\draw[thick, black] (11.8, 0.6) -- (11.8, 1.5);
	 		\draw[thick, black] (12.1, 0.6) -- (12.1, 1.5); 
	 		
	 		\draw[thick, black] (12.4, 0.6) -- (12.4, 1.5);
	 		\draw[thick, black] (12.7, 0.6) -- (12.7, 1.5);
	 		\draw[thick, black] (13, 0.6) -- (13, 1.5);
	 		\draw[thick, black] (13.3, 0.6) -- (13.3, 1.5);
	 		\draw[thick, black] (13.6, 0.9) -- (13.6, 1.5);
	 		\draw[thick, black] (13.9, 0.9) -- (13.9, 1.5);
	 		\draw[thick, black] (14.2, 1.2) -- (14.2, 1.5);
	 		\draw[thick, black] (14.5, 1.2) -- (14.5, 1.5);
	 		
	 		\draw (8.5, 0) node{$\rightarrow$};
	 		\draw (9.5, 0) node{$\psi(\lambda)$=};
	 		\draw (10.15, 0.45) node{$\bullet$}; 
	 		\draw (10.15, 0.15) node{$\bullet$}; 
	 		\draw (10.15, -0.15) node{$\bullet$};
	 		 \draw (10.15, -0.45) node{$\bullet$}; 
	 		 \draw (10.15, -0.75) node{$\bullet$};
	 		  \draw (10.15, -1.05) node{$\bullet$}; 
	 		  \draw (10.15, -1.35) node{$\bullet$}; 
	 		  
	 		  		\draw (10.45, 0.45) node{$\bullet$}; 
	 		  \draw (10.45, 0.15) node{$\bullet$}; 
	 		  \draw (10.45, -0.15) node{$\bullet$};
	 		  \draw (10.45, -0.45) node{$\bullet$}; 
	 		  \draw (10.45, -0.75) node{$\bullet$};
	 		  \draw (10.45, -1.05) node{$\bullet$}; 
	 		  \draw (10.45, -1.35) node{$\bullet$}; 
	 		  
	 		   \draw (10.75, 0.45) node{$\bullet$};
	 		  
	 		 \draw (12.25, 0.75) node{$\bullet$};
	 		\draw (12.55, 0.75) node{$\bullet$};
	 		\draw (12.85, 0.75) node{$\bullet$};
	 		\draw (13.15, 0.75) node{$\bullet$};

	 		 \draw (12.85, 1.05) node{$\bullet$};
	 		 \draw (13.15, 1.05) node{$\bullet$};
	 		 \draw (13.45, 1.05) node{$\bullet$};
	 		 \draw (13.75, 1.05) node{$\bullet$};

	 		 \draw (13.45, 1.35) node{$\bullet$};
	 		 \draw (13.75, 1.35) node{$\bullet$};
	 		 \draw (14.05, 1.35) node{$\bullet$};
	 		 \draw (14.35, 1.35) node{$\bullet$};
	 		 	 		
	 	\end{tikzpicture}
	 \end{figure}

 \newpage
  Case 8: 
 	$s=1,  t=2$ and $\lambda_{1}^{e}-\lambda_{1}^{o} =1.$ 
 	 For any $\lambda\in A_8$ with $n\geq 20$, define 
 	\begin{align*}
 		\psi(\lambda):=(\lambda_1^{e}-3,\lambda_1^{o}-4
 		)	\cup (2^{(\lambda_2^{o}+7)/2}).
 	\end{align*}
 	Note that   if  $n \geq 20 $, then $\lambda_1^{o}\geq 7$.
 	 Hence, 
 	\begin{align*}
 		B_{8}(n)=\{ & \mu \in B_{eu}^{od}(n) \ | v=2, u \geq  4, \mu_1^{o}-\mu_2^{o}=2,
 		\mu_1^{e}=2, \mu_2^{o}-	2u+11>0,\ \mu_1^{o} \geq 5
 		\}.
 	\end{align*}
   In addition, for any $\mu\in B_8(n)$ with $n\geq 20$, 
 	\begin{align*}
 		\psi^{-1}(\mu)=( \mu_1^{o}+3, \mu_2^{o}+4,
 		2f_2 (\mu)-7).
 	\end{align*}

 	Example.  Let $n=22$ and $\lambda=(10,9,3)$. Then $\psi(\lambda)=(10-3,9-4)\cup (2^{5})=(7,5,2,2,2,2,2).$
 	
 	 
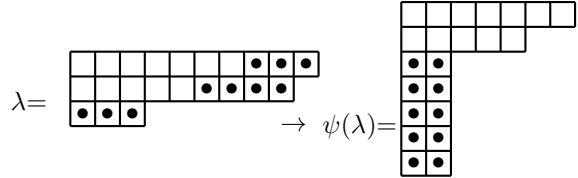
\begin{figure}[htbp] \centering
 		\begin{tikzpicture}[scale=1.1]
 			\draw[thick, black] (4, 0) -- (6.7, 0);
 			\draw[thick, black] (4, 0.3) -- (7, 0.3);
 			\draw[thick, black] (4, 0.6) -- (7, 0.6);
 			\draw[thick, black] (4, -0.3) -- (4.9, -0.3);
 			 
 			\draw[thick, black] (4, -0.3) -- (4, 0.6);
 			\draw[thick, black] (4.3, -0.3) -- (4.3, 0.6);
 			\draw[thick, black] (4.6, -0.3) -- (4.6, 0.6);
 			\draw[thick, black] (4.9, -0.3) -- (4.9, 0.6);
 			\draw[thick, black] (5.2, 0) -- (5.2, 0.6);
 			\draw[thick, black] (5.5, 0) -- (5.5, 0.6);
 			\draw[thick, black] (5.8, 0) -- (5.8, 0.6);
 			\draw[thick, black] (6.1, 0) -- (6.1, 0.6);
 			\draw[thick, black] (6.4, 0)-- (6.4, 0.6);
 			\draw[thick, black] (6.7, 0) -- (6.7, 0.6);
 			\draw[thick, black] (7, 0.3)--(7,0.6);
 			
 			\draw (3.5, 0) node{$\lambda$=};
 			\draw (4.15, -0.15)node{$\bullet$};
 			\draw (4.45, -0.15)node{$\bullet$};
 			\draw (4.75,-0.15)node{$\bullet$};
 			
 			\draw (5.65,0.15) node{$\bullet$};
 			\draw (5.95, 0.15) node{$\bullet$};
 			\draw (6.25, 0.15) node{$\bullet$};
 			\draw (6.55, 0.15) node{$\bullet$};
 			
 			\draw (6.25, 0.45) node{$\bullet$};
 			\draw (6.55, 0.45) node{$\bullet$};
 			\draw (6.85, 0.45) node{$\bullet$};

 			\draw[thick, black] (8, 0) -- (8.6, 0);
 			\draw[thick, black] (8, 0.3) -- (8.6, 0.3);
 			\draw[thick, black] (8, 0.6) -- (9.5, 0.6);
 			\draw[thick, black] (8, 0.9) -- (10.1, 0.9);
 			\draw[thick, black] (8, 1.2) -- (10.1, 1.2);
 			\draw[thick, black] (8, -0.3) -- (8.6, -0.3);
 			\draw[thick, black] (8, -0.6) -- (8.6, -0.6);
 			\draw[thick, black] (8, -0.9) -- (8.6, -0.9);

 			\draw[thick, black] (8, -0.9) -- (8, 1.2);
 			\draw[thick, black] (8.3, -0.9) -- (8.3, 1.2);

 			\draw[thick, black] (8.6, -0.9) -- (8.6, 1.2);
 			\draw[thick, black] (8.9, 0.6) -- (8.9, 1.2);
 			\draw[thick, black] (9.2, 0.6) -- (9.2, 1.2);
 			\draw[thick, black] (9.5, 0.6) -- (9.5, 1.2);
 			\draw[thick, black] (9.8, 0.9) -- (9.8, 1.2);
 			\draw[thick, black] (10.1, 0.9) -- (10.1, 1.2); 
 			
 			 			\draw (6.7,-0.3)node{$\rightarrow$};
 			\draw (7.5, -0.3) node{$\psi(\lambda)$=};
 			\draw (8.15, 0.15) node{$\bullet$};
 			\draw (8.45, 0.15) node{$\bullet$};
 			
 			\draw (8.15, -0.15) node{$\bullet$};
 			\draw (8.45, -0.15) node{$\bullet$};
 			
 			\draw (8.15, -0.45) node{$\bullet$};
 			\draw (8.45, -0.45) node{$\bullet$};
 			
 			\draw (8.15, -0.75) node{$\bullet$};
 			\draw (8.45, -0.75) node{$\bullet$};
 			
\draw (8.15, 0.45) node{$\bullet$};
\draw (8.45, 0.45) node{$\bullet$}; 
 			
 		\end{tikzpicture}
 	\end{figure}

 	Case 9:  $s=1,t= 1$  and $\lambda_{1}^{e}-\lambda_{1}^{o} =1.$ 
 	 Note that for this case,
 	\[
 A_9(4k)=A_9(4k-3)=A_9(4k-2)=\emptyset\qquad  (k\geq 1)
 	\]
 	and 
 		\[
 	A_9(4k-1)=\{(2k,2k-1)\} \qquad  (k\geq 1).
 	\]
 	  For $k\geq 6$,   define 
 	\begin{align*}
 		\psi((2k,2k-1)):=(4k-15,5,3,2,2,2).
 	\end{align*}
 Thus,	
 	\begin{align*}
 		B_{9}(4k-1)=\{(4k-15,5,3,2,2,2)
 		\}
 	\end{align*}
 	and
 		\[
 	B_9(4k)=B_9(4k-3)=B_9(4k-2)=\emptyset\qquad  (k\geq 6).
 	\]
  	Obviously, 
 	\begin{align*}
 		\psi^{-1}((4k-15,5,3,2,2,2))=(2k,2k-1).
 	\end{align*}
 
 Case 10:
 	$s=2$ and $t=1$. 
 For any $\lambda\in A_{10}(n)$
 	 with $n\geq83$,   
 	 define
 	\begin{align*}
 		\psi(\lambda):=(\lambda_2^{e}+7,\lambda_1^{o}+6,5
 		)\cup (2^{(\lambda_1^{e}-18)/2}).
 	\end{align*}
Note  that $\lambda_1^{e}\geq 28$ since $n\geq 83$. Thus, 
 	\begin{align*}
 		B_{10}(n)=\{ & \mu \in B_{eu}^{od}(n) \ | v=3, u\geq 5, \mu_{3}^{o}=5,    \mu_1^e=2,\ 2u+25\geq \mu_1^{o}
 		\}.
 	\end{align*} 
 	In addition, for any $\mu \in B_{10}(n)$ with $n\geq 83$, 
 \begin{align*}
 	\psi^{-1}(\lambda):=(2u+18,\mu_1^{o}-7,\mu_2^{o}-6).
 \end{align*}
  
   Example. Let $n=85$ and $\lambda=(30,28,27)$.
    Then $\psi(\lambda)=(35,33,5,2,2,2,2,2,2)$.
 
 Case 11: $s\geq 3, \ t=1$ and 
 	$   \lambda_1^o=1.$
 	 For any $\lambda \in A_{11}(n)$ with $n\geq 7$, define
   	\begin{align*}
   	\psi(\lambda):=(\lambda_1^{e}+1,\lambda_2^{e},\ldots,\lambda_s^{e}).
   \end{align*}
    Therefore, 
    \begin{align*}
    	B_{11}(n)=\{  \mu \in B_{eu}^{od}(n) \  | u\geq 2,\ v=1 \}.
    \end{align*}
    For any $\mu \in B_{11}(n)$ with $n\geq 7$,
    \[
    	\psi^{-1}(\mu):=(\mu_1^{o}-1,\mu_1^{e},...,\mu_u^{e}, 1).
  \]

  	Example.
  	  Let $n=25$ and $\lambda=(10,8,6,1)$. Then $\psi(\lambda)=(10+1,8,6)=(11,8,6).$

 Case 12:
 	$s=3,\ t =1$ and $ \lambda_1^{o}\geq 3$. 
 For any $\lambda\in A_{12}(n)$ with $n\geq 95$, define  
 	\begin{align*}
 		\psi(\lambda):=(\lambda_2^{e}+7,\lambda_3^{e}+5,\lambda_1^{o}+4)\cup (2^{(\lambda_1^{e}-16)/2})
 	\end{align*}
 	Note that  that if  $n\geq95$, then $\lambda_1^e\geq 24$. Thus, 
 	\begin{align*}
 		B_{12}(n)=\{ & \mu \in B_{eu}^{od}(n) \ |  v=3, \mu_3^o\geq 7, u\geq 4, \mu_{1}^e =2
 	,\ 2u+23\geq \mu_1^{o}	\}.
 	\end{align*}
 	For any $\mu\in B_{12}(n)$ with $n\geq 95$, 
 	 \begin{align*}
 		\psi^{-1}(\mu):=(2u+16,\mu_1^{o}-7,\mu_2^{o}-5,\mu_3^{o}-4).
 	\end{align*}
 
 Example. Let $n=103$ and $\lambda=(26,26,26,25)$. Then $\psi(\lambda)=(26+7,26+5,25+4)\cup(2^{5})=(33,31,29,2,2,2,2,2).$

 Case 13: $s=4, \ t =1$ and $\lambda_1^{o}\geq 3$. 
 For any $\lambda\in A_{13}(n)$ with $n\geq 159$, 	define 
 	\begin{align*}
 		\psi(\lambda):=(\lambda_2^{e}+7,\lambda_3^{e}+5,\lambda_4^{e}+3,\lambda_1^{o}+2,3)\cup (2^{(\lambda_1^{e}-20)/2}).
 	\end{align*}
 	Note that 	if $n\geq159$, then $\lambda_1^e\geq 32$.
 	Therefore, 
 	\begin{align*}
 		B_{13}(n)=\{ & \mu \in B_{eu}^{od}(n) \ | u\geq 6, \ v=5, \ \mu_5^o=3, \  \mu_{1}^e=2
 		, 2u+27\geq \mu_{1}^{o}\}.
 	\end{align*}
 	For any $\mu \in B_{13}(n)$ with $n\geq 159$,
 	\begin{align*}
 		\psi^{-1}(\mu)=(2u+20,\mu_1^{o}-7,\mu_2^{o}-5,\mu_3^{o}-3,\mu_4^{o}-2).
 	\end{align*}
 	
 	Example. Let $n=179$ and $\lambda=(38,36,36,36,33)$. 
 	 Then $\psi(\lambda)=(36+7,36+5,36+3,33+2,3)\cup(2^{9})=(43,41,39,35,3,2,2,2,2,2,2,2,2,2).$

  Case 14:  $s=5,\ t =1$ and $ \lambda_1^{o}\geq 3$.
   	For any $\lambda \in A_{14}(n)$ with $n\geq 227$,  define 
  	\begin{align*}
  		\psi(\lambda):=(\lambda_2^{e}+9,\lambda_3^{e}+7,\lambda_4^{e}+5,\lambda_4^{e}+3,\lambda_1^{o}+2)\cup (2^{(\lambda_1^{e}-26)/2}).
  	\end{align*}
  Note that if 	$n\geq227$,  then $\lambda_1^e\geq 38$.
  	Thus, 
  	\begin{align*}
  		B_{14}(n)=\{ & \mu \in B_{eu}^{od}(n) \ | u\geq 6,\  
  		 v=5,\  \mu_5^o\geq 5, \  \mu_{1}^e =2, 2u+35\geq \mu_{1}^{o}
  		\}.
  	\end{align*}
  	For any $\mu \in B_{14}(n)$ with $n\geq 227$, 
  	\begin{align*}
  		\psi^{-1}(\mu)=(2u+26,\mu_1^{o}-9,\mu_2^{o}-7,
  		\mu_3^{o}-5,\mu_4^{o}-3,\mu_5^{o}-2).
  	\end{align*}

 Case 15: 
  	$6\leq s\leq 10,\ t =1$ and $ \lambda_1^{o}\geq 3$. 
  	For any $\lambda\in A_{15}(n)$ with $n\geq 373$, define 
  	\begin{align*}
  		\psi(\lambda):=(\lambda_2^{e}+5,\lambda_3^{e}+3,\lambda_4^{e}+1,\lambda_5^{e},\lambda_6^{e},...,\lambda_s^{e},
  		\lambda_1^{o}+1
  		)\cup (2^{(\lambda_1^{e}-10)/2})
  	\end{align*}
  Note that if 	  $n\geq 373$,  then  $\lambda_1^e \geq  34$. Therefore, 
  	\begin{align*}
  		B_{15}(n)=\{ & \mu \in B_{eu}^{od}(n) \ | f_2(\mu)>12,\  v=3,   \  	
  		 \mu_{u-f_2(\mu)}^{e}\geq 4, \  3 \leq 	u-f_2(\mu)\leq  7, \ 2f_2(\mu)+15\geq \mu_1^{o} 
  		\}.
  	\end{align*}
  	For any $\mu\in B_{15}(n)$
  	 with $n\geq 373$, 
  	\begin{align*}
  		\psi^{-1}(\mu):=(2f_2(\mu)+10,\mu_1^{o}-5,\mu_2^{o}-3,\mu_3^{o}-1,\mu_1^{e},\mu_2^{e},\ldots,	\mu_{u-f_2(\mu)-1}^{e},
  		\mu_{u-f_2(\mu)}^{e}-1).
  	\end{align*}

 Case 16: 
  	$s\geq 11,\ t=1,\  \lambda_1^{o}\geq 3$ and $ \lambda_1^{e}-\lambda_2^{e}\leq 10$.
  	 For any $\lambda\in A_{16}(n)$ with $n\geq 47$, define 
  	\begin{align*}
  		\psi(\lambda):=(\lambda_1^{e}+5,\lambda_2^{e}+3,\lambda_3^{e}+1,\lambda_4^{e},\ldots,\lambda_{s-4}^{e},\lambda_{s-3}^{e}-2,\lambda_{s-2}^{e}-2,\lambda_{s-1}^{e}-2,\lambda_s^{e}-2,
  		\lambda_1^{o}-1).
  	\end{align*}
  	Therefore, 
  	\begin{align*}
  		B_{16}(n)=\{ & \mu \in B_{eu}^{od}(n) \ | u\geq 9,\ 
  		 v=3, \  f_2(\mu)\leq 5, \ 	   	\mu_1^{o}-\mu_2^{o}\leq 12
  		,\ \mu_{u-5}^{e}-\mu_{u-4}^{e}\geq 2\}.
  	\end{align*}
  	For any $\mu\in B_{16}(n)$ with $n\geq 47$, 
  	\begin{align*}
  		\psi^{-1}(\mu)=(\mu_1^{o}-5,\mu_2^{o}-3,\mu_3^{o}-1,\mu_1^{e},\mu_2^{e},...\mu_{u-5}^{e},\mu_{u-4}^{e}+2,\mu_{u-3}^{e}+2,\mu_{u-2}^{e}+2,	\mu_{u-1}^{e}+2,
  		\mu_{u}^{e}+1
  		).
  	\end{align*}

 Case 17:
 	$ s\geq11, \ t=1,\ \lambda_1^{o}\geq 3$ and $\lambda_1^{e}-\lambda_2^{e}\geq 12$. 
  For any $\lambda \in A_{17}(n)$ with $n\geq 59$,
    	define  
 	\begin{align*}
 		\psi(\lambda):=(\lambda_1^{e}-7,\lambda_2^{e}+3,\lambda_3^{e}+1,\lambda_4^{e},...,\lambda_{s-4}^{e},\lambda_{s-3}^{e}-2,\lambda_{s-2}^{e}-2,\lambda_{s-1}^{e}-2,\lambda_s^{e}-2,
 		\lambda_1^{o}-1,2,2,2,2,2,2).
 	\end{align*}
 	Therefore, 
 	\begin{align*}
 		B_{17}(n)=\{ & \mu \in P_{eu}^{od}(n) \ | u\geq 15,\ v=3, \  6 \leq  f_2(\mu)\leq 11, \mu_{u-11}^{e}-\mu_{u-10}^{e}\geq 2
 		\}.
 	\end{align*}
 	For any $\mu \in B_{17}(n)$
 	 with $n\geq 59$,
 	\begin{align*}
 		\psi^{-1}(\mu):=(\mu_1^{o}+7,\mu_2^{o}-3,\mu_3^{o}-1,\mu_1^{e},\mu_2^{e},\ldots,\mu_{u-11}^{e},\mu_{u-10}^{e}+2,\mu_{u-9}^{e}+2,\mu_{u-8}^{e}+2,	\mu_{u-7}^{e}+2,
 		\mu_{u-6}^{e}+1).
 	\end{align*}
 	
Therefore, $\psi$ is a bijection from $A_j(n)$ to $B_j(n)$ when $n\geq 373$.
 Thus for $n\geq 373$,
 \begin{align}
 	|A_j(n)|=|B_j(n)|. \label{3-6}
 \end{align}
 Theorem \ref{Th-1} follows from \eqref{3-4}--\eqref{3-6}. This completes
  the proof. \qed

\section{ Concluding remarks}

As seen in Introduction, partitions with parts separated by parity
  have received
considerable attention in recent years.
 Recently, Bringmann, Craig and Nazaroglu \cite{Bringmann} 
 studied the asymptotic behavior of the eight partition functions  proved  several inequalities  for sufficiently large $n$.   They  also  asked for combinatorial 
 proofs of those inequalities at the end of their paper.  
 In this paper, we show  that  an inequality   on 
 partitions separated by parity  holds for $n\geq 373$ by a    combinatorial method. This  settles an open problem  
 posed  by Bringmann, Craig and Nazaroglu. 
  Unfortunately,
   we can not give combinatorial proofs of \eqref{1-2}
    and \eqref{1-3}.  The combinatorial proofs of \eqref{1-2}
    and \eqref{1-3} will likely
    require a different approach.
    
    In \cite{Ballantine},  Ballantine 
     and Welch  provided
     a partition-theoretic interpretation
     of  $p_{od}^{ed}(n)-p_{ed}^{od}(n)$. Therefore,
      it would be interesting to present   
      a partition-theoretic interpretation
      of $ p_{eu}^{od}(n)-p_{od}^{eu}(n)$ for $n\geq 50$.

\section*{Statements and Declarations}

\noindent{\bf Acknowledgments.}  The authors  would like to
express their sincere gratitude to
Shishuo Fu  for his helpful suggestions.   This work was supported by
the National Science Foundation of
China  (Grant Nos. 12401434 and 
12371334) and the Natural Science Foundation of
Jiangsu Province of China (Grant no.
BK20221383).

\noindent{\bf  Author Contributions.}
  The  authors have equally contributed to the manuscript.

\noindent{\bf Declaration of competing interest.}
The authors declare that they have  no
known competing financial
interests or personal relationships
that could have appeared to
influence the work reported in this paper.

\noindent{\bf Data availability.}
No data was used for the research
described in the article.

\end{document}